\pdfoutput=1
\documentclass[lettersize]{article}
\usepackage{sampta,amsmath,amsthm,epsfig,fancyhdr}

\usepackage{url}
\def\pathPic{Illustrations/}

\def\cE{\mathcal E}
\def\cH{{\cal H}}
\def\cM{\mathcal M}
\def\cT{{\cal T}}
\def\P{{\rm \hbox{I\kern-.2em\hbox{P}}}}
\def\H{{\rm \hbox{I\kern-.2em\hbox{H}}}}
\def\R{{\rm \hbox{I\kern-.2em\hbox{R}}}}
\newcommand{\be}{\begin{equation}}
\newcommand{\ee}{\end{equation}}
\newcommand\ssep{; \,}% small separation.
\DeclareMathOperator\interp{{\rm I}}
\def\ve{\varepsilon}
\newcommand{\iref}[1]{(\ref{#1})}
\newcommand\trans{\mathrm T}
\DeclareMathOperator\disc{disc}
\DeclareMathOperator\Id{Id}
\def\proof{{\noindent \bf Proof: }}
\def\sq{\hfill $\diamond$}
\def\TEq{{T_{\rm eq}}}
\def\cPi{{\text{\large$\boldsymbol\pi$}}}

\newtheorem*{theorem*}{Theorem}

\begin{document}\sloppy

% Example definitions.
%\def\x{{\mathbf x}}
%\def\L{{\cal L}}
%\def\sampta{SampTA~}

% Title, name and affiliations
%\title{Piecewise quadratic anisotropic finite elements:\\ optimal triangle or optimal aspect ratio?}
\title{The optimal aspect ratio for piecewise quadratic\\ anisotropic finite element approximation} % (Extended Abstract)}

\name{Jean-Marie Mirebeau}
\address{Laboratoire Jacques Louis Lions, Universit\'e Pierre et Marie Curie, Paris, France\\
E-mail: mirebeau@ann.jussieu.fr}

\maketitle

\begin{abstract}
Mesh adaptation for finite element approximation is a procedure used in numerous applications. The use of thin and long \emph{anisotropic} triangles improves the efficiency of the procedure. 

When piecewise linear finite elements are used, the aspect ratio for mesh adaptation is generally dictated by the absolute value of the (estimated) hessian matrix of the approximated function. We give in this paper the corresponding aspect ratio for piecewise quadratic finite elements. 
\end{abstract}


\section{Introduction}
\label{sec:intro}
Consider a bounded polygonal domain $\Omega\subset \R^2$, a sufficiently smooth function $f:\Omega\to \R$, and an integer $m\geq 2$. We introduce the problem of \emph{optimal mesh adaptation} 
\be
\label{eqOptF}
\min \{\#(\cT) \ssep \cT \text{ s.t. }\| \nabla (f-\interp_\cT^{m-1} f)\|_{L^2(\Omega)}\leq \ve\},
\ee
where $\cT$ stands for an arbitrary triangulation of $\Omega$, and $\#(\cT)$ for its cardinality. % andand the optimization is among all triangulations $\cT$ of $\Omega$, and $\#(\cT)$ denotes the cardinality of $\cT$. 
Here $\interp_\cT^{m-1}$ denotes the Lagrange interpolation operator onto finite elements of degree $m-1$ on $\cT$.

In practical applications, the problem \iref{eqOptF} is 
%hardly tractable 
generally intractable 
for at least three reasons. 1: The function $f$ may have complicated local features, difficult to analyze. We thus first make a local analysis based on Taylor developments. 2: The collection of triangular meshes of $\Omega$ is a combinatorial set and problems such as \iref{eqOptF} are typically NP-complete (after discretization). We avoid this problem by first considering the case of a single triangle. 3: Currently available anisotropic mesh generation algorithms only give control on the aspect ratio and orientation of the generated triangles, but not on their other features. We thus only optimize this aspect ratio.

\section{An optimization problem}

%We first focus on the local behavior of $f$, and for that purpose we define the degree $m:=k+1\geq 2$.
We denote by $\P_{m-1}$ the space of bivariate polynomials of degree $\leq m-1$, and by $\H_m$ the space of homogeneous  polynomials of degree $m$.
If $f\in C^m(\Omega)$, if $z\in \Omega$ is fixed and if $h\in \R^2$ is small, then locally
\be
\label{eqTaylor}
f(z+h) = \mu_z(h)+ \pi_z(h) +o(|h|^m),
\ee
for some $\mu_z\in \P_{m-1}$ and $\pi_z\in \H_m$. 
If $T$ is a sufficiently small triangle, we thus have at least heuristically on $T$
\be
\label{eqLocF}
\nabla (f-\interp_T^{m-1}f) \simeq \nabla (\pi_z - \interp_T^{m-1} \pi_z),
\ee
since the Lagrange interpolation operator $\interp_T^{m-1}$ on the triangle $T$ reproduces the elements of $\P_{m-1}$.

%We now introduce 
For any triangle $T$ and any $f\in H^1(T) \cap C^0(T)$, we define the averaged $H^1$ interpolation error $e_T(f)_m$ as follows
%
%We define $m:=k+1\geq 2$, and for each triangle $T$ and function $f\in C^0(T)$ we define the averaged $H^1$ interpolation error on $T$ as follows
$$
e_T(f)_m^2 := \frac 1 {|T|} \int_T |\nabla (f-\interp_T^{m-1} f)|^2.
$$
The local counterpart of \iref{eqOptF} is the problem of the \emph{optimal triangle} : find for all $\pi\in \H_m$
%We introduce a ``local'' counterpart of \iref{eqOptF} as follows : for all  $\pi\in \H_m$
% In view of \iref{eqLocF}, the counterpart of \iref{eqOptF} based on local behavior and for a single triangle is the following optimization problem, which is posed for all
\be
\label{optT}
\sup\{ |T| \ssep T \text{ s.t. } e_T(\pi)_m \leq 1\}.
\ee
%where the supremum is taken among all bi-dimensional triangles $T$. 
Indeed the cardinality of a triangulation is inversely proportional to the area of its elements. This approach is developed in Chapter 2 of \cite{thesisMi} and leads to asymptotically optimal error estimates of \iref{eqOptF} as $\ve \to 0$ (or more precisely estimates of $\ve$ as $\#(\cT) \to \infty$, which is equivalent). Unfortunately these estimates are not completely realistic for applications because currently available numerical anisotropic mesh generators only control the \emph{aspect ratio and orientation} of the generated triangles.

For each triangle $T$, of vertices $v_1$, $v_2$ and $v_3$, we denote by $z_T := (v_1+v_2+v_3)/3$ its barycenter. We denote by $S_2^+$ the collection of $2\times 2$ symmetric positive definite matrices, and we define a matrix $\cH_T\in S_2^+$ by the equality 
$$
\cH_T^{-1} := \frac 2 3 \sum_{1\leq i \leq 3} (v_i-z_T) (v_i-z_T)^\trans.
$$
If $A$ is an invertible $2\times 2$ matrix and if $T'$ is mapped onto $T$ by the linear map $z\mapsto Az$, then one easily checks that 
\be
\label{eqInvHT}
\cH_{T'} = A^\trans \cH_T A.
\ee
By construction the triangle $\TEq$ of vertices $(\cos(2 k  \cPi/3), \sin (2 k \cPi/3))_{0\leq k \leq 2}$ satisfies $\cH_\TEq = \Id$. Combining these two properties, Proposition 5.1.3 in \cite{thesisMi} establishes that for any triangle $T$
%If $T'$ is the image of $T$ by a linear map $z\mapsto Az$  easily checka
%Several properties of this matrix are established in \cite{thesisMi}, Proposition 5.1.3.   that 
$$
|T|\sqrt{\det \cH_T} = |\TEq|, %3\sqrt 3/4,
$$
and that there exists a rotation $U$ (depending on $T$) such that 
\be
\label{eqMapTTeq}
z\mapsto U \cH_T^\frac 1 2 (z-z_T)
\ee
maps $T$ onto $\TEq$ (the power $\alpha$ of a symmetric positive definite matrix is obtained by elevating the eigenvalues to the power $\alpha$ in a diagonalization). 
%where $\TEq$ is the triangle of vertices $(\cos(2 k  \cPi/3), \sin (2 k \cPi/3))_{0\leq k \leq 2}$.
Furthermore 
the ellipse of minimal volume containing $T$ is $\cE_T := \{z \ssep (z-z_T)^\trans \cH_T (z-z_T)\}$, see Fig 1. The matrix $\cH_T$ thus encodes the area, the aspect ratio and the orientation of $T$.
\begin{figure}
\begin{tabular}{ccc}
\includegraphics[width = 2cm,height=2cm]{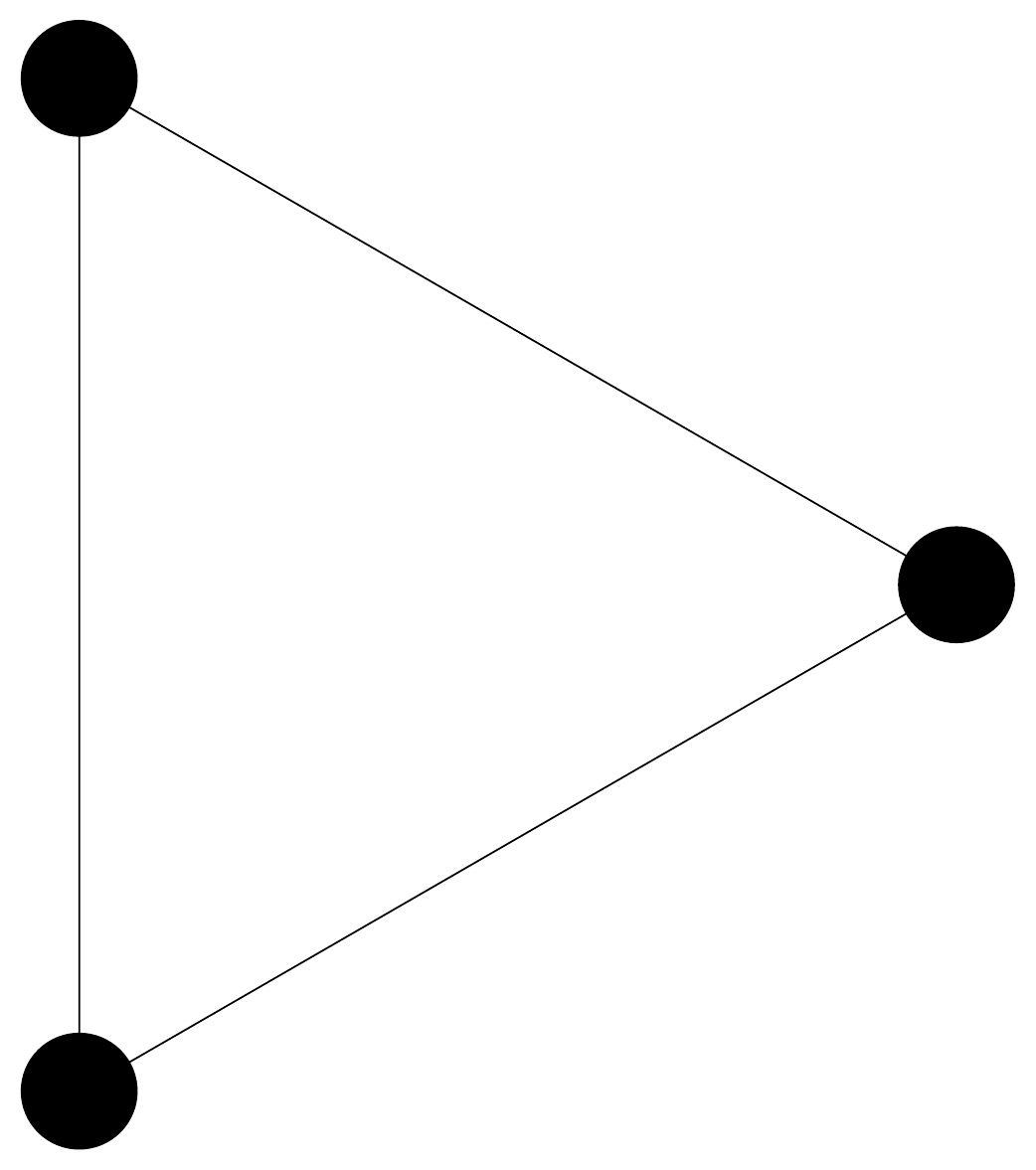}
& \includegraphics[width = 2cm,height=2cm]{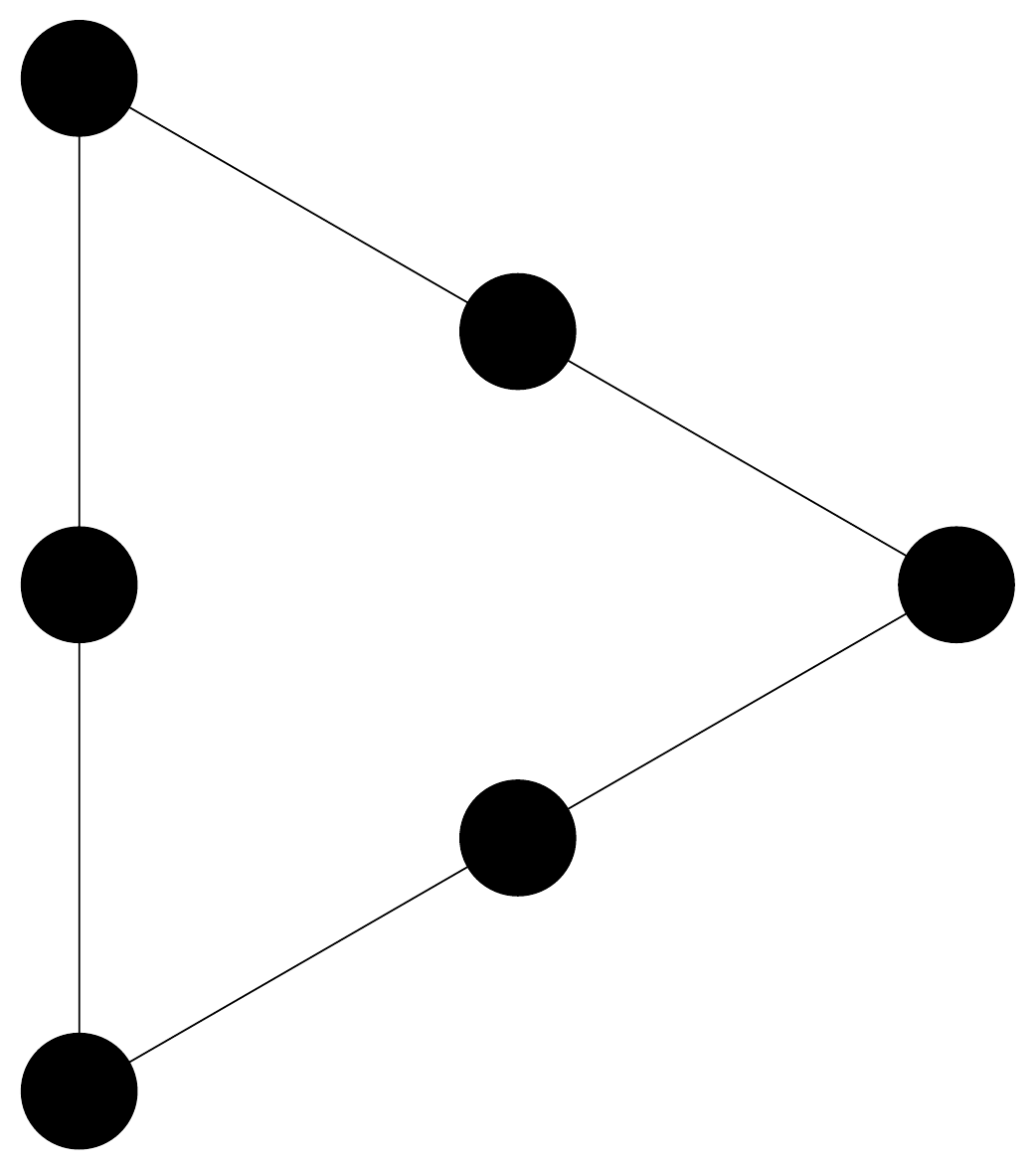}
& 
{\raise 8mm \hbox{
\begin{tabular}{c}
%\vspace{15mm}
\includegraphics[width=3cm,height=1.2cm]{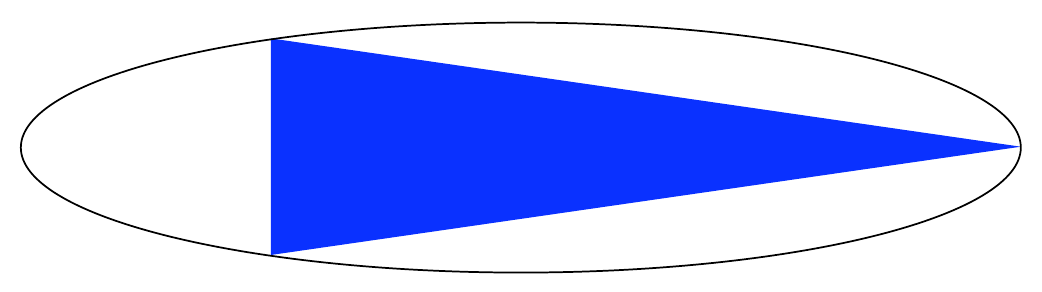}\\
%\vspace{15mm}
\includegraphics[width=3cm,height=1.2cm]{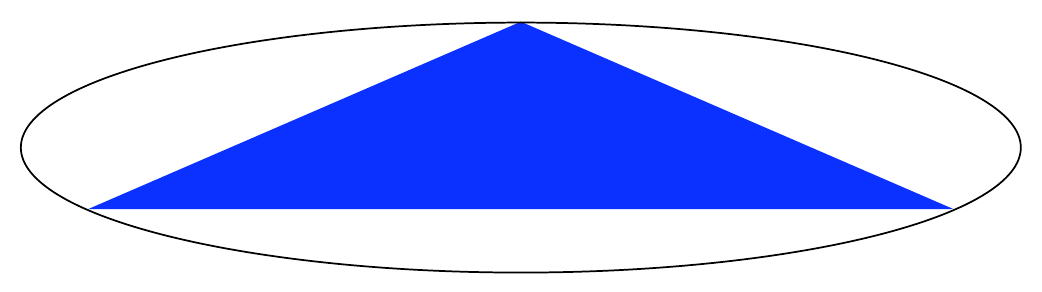}
\end{tabular}
}
}
\end{tabular}
\caption{Lagrange interpolation points for $\P_1$ and $\P_2$ finite elements (left), triangle $T$ and associated ellipse $\cE_T$ (right).}
\end{figure}

For each $M\in S_2^+$ and each $\pi\in \H_m$ we define %introduce the quantity 
$$
e_M(\pi)_m := \sup\{ e_T(\pi)_m \ssep T\text{ s.t. } \cH_T = M\}.
$$
We finally introduce for each $\pi\in \H_m$ the problem of the \emph{optimal aspect ratio} for $\P_{m-1}$ interpolation
\be
\label{optM}
\inf\{\det M \ssep M \in S_2^+ \text{ s.t. }e_M(\pi)_m \leq 1\}.
\ee

\section{Main result}

Our main result is the solution of the optimization problem \iref{optM} in the case of piecewise linear and piecewise quadratic finite elements. The piecewise quadratic case is entirely new and gives a well founded answer to a long standing question: which aspect ratio, depending on the third derivatives of the approximated function, should be used in finite element software that combine anisotropy and $\P_2$ elements ?

We first introduce some notation.
We equip the vector space $\H_m$ with the norm
$$
\|\pi\|:= \sup_{|u|\leq 1} |\pi(u)|.
$$
For each $\pi \in \H_2$, $\pi = a x^2+2bxy+cy^2$, we define 
$$
[\pi] = \left(
\begin{array}{cc}
a & b\\
b & c
\end{array}
\right).
$$
The absolute value of a symmetric matrix (resp. the square root of a non negative symmetric matrix) is obtained by taking the absolute value (resp. square root) of the eigenvalues in a diagonalization.
For each $\pi\in \H_2$ we set 
$$
\cM_2(\pi) := \|\pi\|  \ |[\pi]|.
$$
For each $\pi\in \H_3$, $\pi = a x^3+ 3 b x^2y+ 3 c x y^2+ d y^3$, we set 
$$
\cM_3(\pi) := \sqrt{ [\partial_x \pi]^2+[\partial_y\pi]^2}+ \left(\frac{-\disc \pi}{\|\pi\|}\right)_+^{\frac 1 3} \Id,
$$
where $\disc \pi := 4(ac-b^2)(bd-c^2) - (ad-bc)^2$ and $\lambda_+ := \max\{\lambda,0\}$.

\begin{theorem*}
For $m\in \{2,3\}$ the map $\pi\in \H_m \to \cM_m(\pi)$ is a \emph{near-minimizer} of the problem \iref{optM} in the following sense. If $\pi$ is non-univariate then $\cM_m(\pi)$ is non-degenerate. Furthermore there exists a constant $C$, independent of $\pi$, such that $e_{\cM_m(\pi)}(\pi)_m \leq C$ and 
$$
\det \cM_m(\pi) \leq C \inf\{\det M \ssep M\in S_2^+ \text{ s.t. }  e_M(\pi)_m \leq 1\}.
$$
\end{theorem*}

\proof 
The integer $m\in\{2,3\}$ is fixed, and we denote for each $\pi\in \H_m$ 
$$
\|\nabla \pi\| := \sup_{|u|\leq 1} |\nabla \pi(u)|.
$$
For %each $\pi\in \H_m$ and 
each $2\times 2$ matrix $A$ we denote by $\pi \circ A$ the element of $\H_m$ defined by
$
(\pi \circ A)(u) := \pi(A(u)), \ u\in \R^2.
$
We recall that 
$$
\nabla (\pi \circ A) (u)  = A^\trans \nabla \pi (A (u)), \ u\in \R^2,
$$
which implies for any rotation $U$
\be
\label{eqChgNorm}
\|\nabla \pi\| = \|\nabla (\pi\circ U)\|.
\ee
%linear 
%We denote for 
The main difficulty of this proof is to show that there exists a constant $C\geq 1$ such that for all $\pi \in \H_m$ and all $M\in S_2^+$ one has 
\be
\label{eqErMNorm}
C^{-1} e_M(\pi)_m \leq \|M\|^\frac 1 2 \| \nabla (\pi \circ M^{-\frac 1 2})\| \leq C e_M(\pi)_m.
\ee
Assume that this point is established. %The problem \iref{optM} of optimal aspect ratio is thus equivalent up to the multiplicative constant $C$ to the problem
%\be
%\label{optM2}
%\inf\{\det M \ssep \|M\|^\frac 1 2 \| \nabla (\pi \circ M^{-\frac 1 2})\|\}.
%\ee
Proposition 6.5.4 in \cite{thesisMi}, states that the map $\pi\mapsto \cM_m(\pi)$ is a near-minimizer for the optimization problem 
$$
\inf\{\det M \ssep M \in S_2^+ \text{ s.t. } \|M\|^\frac 1 2 \| \nabla (\pi \circ M^{-\frac 1 2})\|\},
$$
%\iref{optM2},
in the same sense as in the statement of this theorem. Combining this result with the equivalence  \iref{eqErMNorm}, and using the homogeneity of $\pi$, we immediately conclude the proof of this theorem.
% which immediately concludes the proof using the equivalence \iref{eqErMNorm}.

We thus turn to the proof of \iref{eqErMNorm}.
Our first observation is that there exists a constant $C_0$ such that for all $\pi\in \H_m$
\be
\label{eqLow}
e_\TEq (\pi)_m := \sqrt{\frac {1} {|\TEq|} \int_\TEq |\nabla (\pi-\interp_\TEq^{m-1} \pi)|^2} \leq C_0\|\nabla \pi\|,
\ee
indeed the left and right hand side are norms on $\H_m$.

Consider a symmetric matrix $M\in S_2^+$ and a triangle $T$ such that $\cH_T = M$. According to \iref{eqMapTTeq} there exists a rotation $U$ such that the image of $T$ by the map 
$
z\mapsto U M^\frac 1 2 (z-z_T) 
$
is the triangle $\TEq$. Injecting this change of variables in \iref{eqLow} we obtain %Using this change of variables we obtain %and the fact that the interpolation operator is invariant under we obtain 
$$
\sqrt{\frac {1} {|T|} \int_T |U M^{-\frac 1 2}\nabla (\pi-\interp_T^{m-1} \pi)|^2} \leq C_0 \|\nabla (\pi \circ (M^{-\frac 1 2} U^{-1}))\|.
$$
Observing that $\|Av\| \geq \|A^{-1}\|^{-1} |v|$ for any invertible $2\times 2$ matrix $A$ and vector $v\in \R^2$, and recalling \iref{eqChgNorm},
%, and that 
%\be
%\label{eqChgNorm}
%\|\nabla (\pi' \circ U')\| = \| \nabla \pi'\|
%\ee
%for any $\pi' \in \H_m$ and any rotation $U'$, 
we obtain
$$
\|M\|^{-\frac 1 2} e_T(\pi)_m \leq C_0 \| \nabla (\pi \circ M^{-\frac 1 2})\|.
$$
Taking the supremum of the left hand side among all triangles $T$ such that $\cH_T = M$ we establish the left part of \iref{eqErMNorm}, provided that $C\geq C_0$.

We now remark that there exists a constant $C_1$ such that for all $\pi \in \H_m$
\be
\label{eqUp}
\|\nabla \pi\| \leq C_1 \sqrt{\frac {1} {|\TEq|} \int_\TEq |\partial_x (\pi-\interp_\TEq^{m-1} \pi)|^2}.
\ee
Indeed assume that the right hand side vanishes.
%$\|\nabla (\pi \circ (M^{-\frac 1 2} U^{-1}))\| = \|\nabla (\pi \circ M^{-\frac 1 2})\|$
%which immediately implies that 
%$$
%e_T(\pi) \leq 
%$$
%
%Our first observation is that there exists a constant $C\geq 1$ such that for all $\pi \in \H_m$ % and for that purpose we remark that there exists a constant $C$ such that following quantities are norms on $\H_m$, hence are equivalent
%%$$
%%\|\nabla \pi\|, \  e_\TEq(\pi)_m \text{ and } \frac 1 {\TEq }
%%$$
%$$
%\frac {C^{-2}} {|\TEq|} \int_\TEq |\nabla (f-\interp_\TEq f)|^2 \leq \|\nabla \pi\|^2 \leq \frac {C^2} {|\TEq|} \int_\TEq |\partial_x (\pi - \interp_\TEq^{m-1} \pi)|^2.
%$$
%Indeed the square root of these three quantities is a norm on $\H_m$, which is a vector space of finite dimension. For instance if the right hand size vanishes, 
Then $\mu := \pi-\interp_\TEq^{m-1} \pi$ is a polynomial of degree $m$ depending only on the variable $y$, and which vanishes on the Lagrange interpolation points of $\TEq$, see Fig1. Hence $\mu$ vanishes for $y=\pm \sqrt 3/2$ and $y=0$ if $m=2$ (resp. $y=\pm \sqrt 3/2$, $y = \pm \sqrt 3/4$ and $y=0$ if $m=3$). Therefore $\mu = 0$ which implies that $\pi = 0$. %The Lagrange interpolation points on $\TEq$ are illustrated on Fig 1.
Both sides of \iref{eqUp} are thus equivalent norms on the vector space $\H_m$.

We consider a diagonalization of a symmetric matrix $M^\frac 1 2$, $M\in S^2_+$, \vspace{-2mm}
$$
M^\frac 1 2 = U^\trans 
D U, \ \
D = 
\left(
\begin{array}{cc}
\alpha &0\\
0 & \beta
\end{array}
\right),
$$ 
where $U$ is a rotation and $\alpha = \|M\|^\frac 1 2$. Consider the triangle $T$ which is mapped onto $\TEq$ by the change of coordinates
$$
z\mapsto UM^\frac 1 2 z = D U z,
$$
and thus satisfies $\cH_T = M$ according to \iref{eqInvHT}.
Injecting this change of variables into \iref{eqUp} we obtain 
$$
\|\nabla (\pi\circ (M^{-\frac 1 2} U^{-1})) \| \leq C_1 \sqrt{\frac {1} {|T|} \int_T |\alpha ^{-1} v \cdot \nabla  (\pi-\interp_T^{m-1} \pi)|^2},
$$
where $v:= U^{-1} e_x$,  $e_x := (1,0)$, and where we used for the $\partial_x$ derivative that $U^{-1} D^{-1} e_x= \alpha ^{-1} v$. Recalling that $\alpha = \|M\|^\frac 1 2$, $|v|=1$, and using \iref{eqChgNorm} we obtain 
\begin{eqnarray*}
\|M\|^\frac 1 2 \|\nabla (\pi\circ M^{-\frac 1 2}) \| &\leq& C_1 \sqrt{\frac {1} {|T|} \int_T | \nabla (\pi-\interp_T^{m-1} \pi)|^2}\\
 &=& C_1 e_T(\pi)_m  \leq C_1 e_M(\pi)_m.
\end{eqnarray*}
This concludes the proof of \iref{eqErMNorm} with $C:=\max\{C_0, C_1\}$, hence the proof of this theorem.
%Consider a matrix $M\in S_2^+$ and a triangle $T$ such that $\cH_T = M$. According to \cite{thesisMi}, Proposition 5.1.3 there exists a rotation $U$ such that the image of $T$ by the map 
%$$
%z\mapsto U M^\frac 1 2 (z-z_T) 
%$$
%is the triangle $\TEq$. Using this change of variables we obtain 
%$$
%C^{-2} 
%$$
\sq
%\proof
%Reformulation under the form $\|M\|^\frac 1 2 \| \pi \circ M^{-\frac 1 2}\|$, then thesis.
%\sq

\section{Applications and conclusion}

Consider a function $f$ for which one desires to solve, at least heuristically, the optimization problem \iref{eqOptF}. Assume that some estimate of 
$$
\pi_z := \sum_{k+l = m} \frac{\partial^m f}{\partial^k x \partial^l y}(z) \frac {x^k}{k!}\frac {y^l}{l!}
$$
is known at each point $z\in \Omega$, and define a \emph{riemannian metric} $H$ on $\Omega$ as follows
\be
\label{defH}
H(z) := \lambda (\det \cM_m(\pi_z))^{-\frac {1}{2m}} \cM_m(\pi_z),
\ee
where $\lambda>0$ is a constant
(this expression needs to be slightly modified if $\pi_z$ vanishes or is univariate for some values of $z$, in order to ensure that $H \in C^0(\Omega, S_2^+)$).
Some mesh generators such as \cite{FreeFem} can, at least heuristically, and provided $H$ has sufficient regularity, produce a mesh $\cT$ of $\Omega$ such that $C^{-1} H(z) \leq \cH_T \leq C H(z)$ for each $T\in \cT$ each $z\in T$, where $C$ is a constant not too large. In other words the aspect ratio of the elements of $\cT$ is dictated by the metric $H$. Some rigorous results in this direction can be found in Chapter 5 of \cite{thesisMi}. %[Thesis], and are based on [Shew].

In the expression \iref{defH} the matrix $\cM_m(\pi_z)$ ensures that the elements of $\cT$ have the optimal aspect ratio, while the scalar factor $(\det \cM_m(\pi_z))^{-\frac {1}{2m}}$ guarantees that the interpolation error is equidistributed among the elements of $\cT$ (a general principle in adaptive approximation). %, and the constant $\lambda$ 

We conducted some numerical experiments using \cite{FreeFem}  and for the synthetic function
\be
\label{syntheticF}
f(x,y) := \tanh(10 (\sin(5y)-2x)) + x^2 y+ y^3
\ee
on the domain $\Omega := (-1,1)^2$. They illustrate the improvement offered by anisotropic mesh adaptation, both in the case of $\P_1$ and $\P_2$ elements, for a triangulation of cardinality $500$.
$$
\begin{array}{cccc}
\#(\cT) = 500 & \text{Uniform} & \text{Isotropic} & \text{anisotropic} \\
\|\nabla(f-\interp_\cT^1 f)\|_{L^2} & 110 & 51 & 11\vspace{0.5mm}\\
\|\nabla(f-\interp_\cT^2 f)\|_{L^2} & 79   & 14 & 0.88
\end{array}
$$
Our next objective is to combine our analysis with an adaptive anisotropic mesh refinement procedure, for a partial differential equation solved with $\P_2$ finite elements. The optimization problem \iref{eqOptF} is particularly relevant in the case of elliptic equations. %For elliptic equations the problem \iref{eqOptF} is indeed highly relevant.

\begin{figure}
\centering
\vspace{-2mm}
\hspace{-27mm}
\includegraphics[width=5cm,height=3cm]{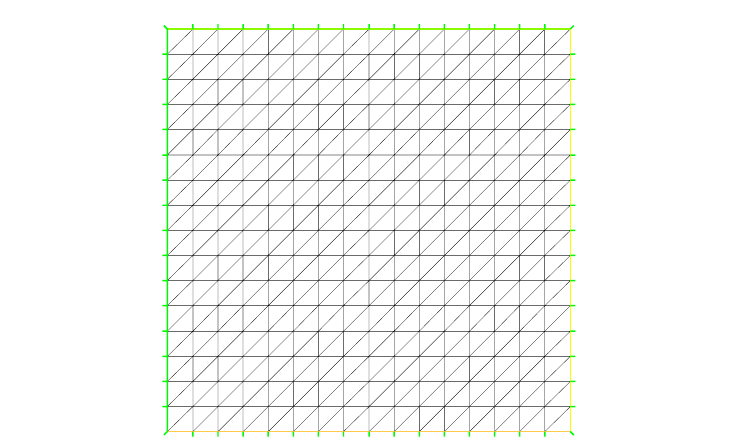}
\hspace{-20mm}
\includegraphics[width=5cm,height=3cm]{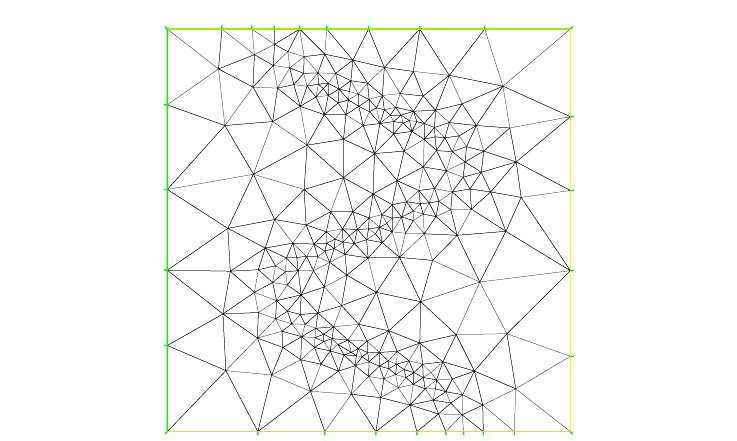}
\hspace{-20mm}
\includegraphics[width=5cm,height=3cm]{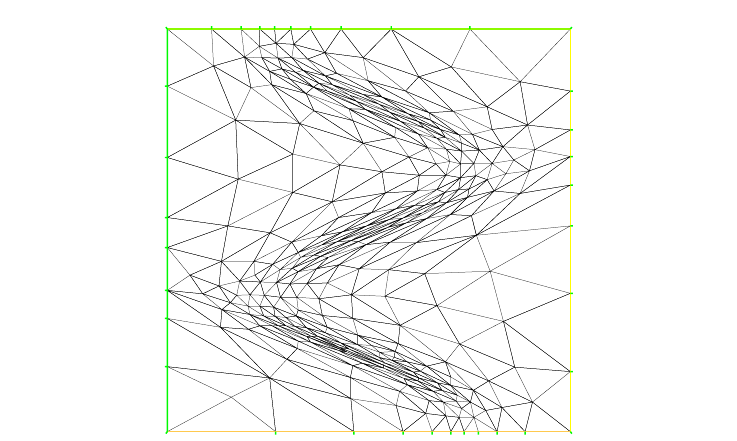}
\hspace{-30mm}

\hspace{-27mm}
\includegraphics[width=5cm,height=3cm]{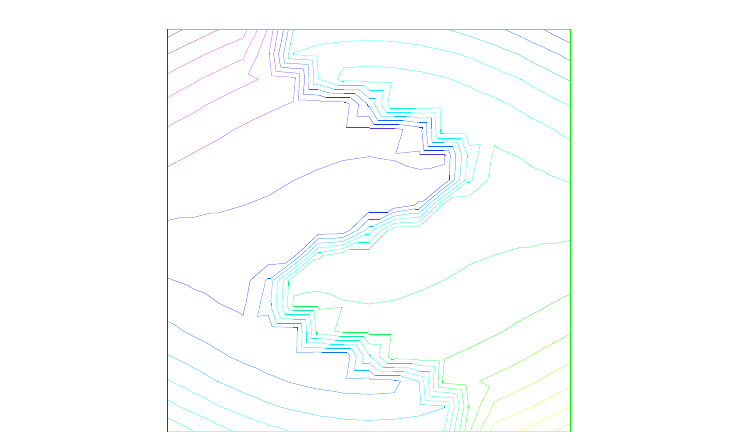}
\hspace{-20mm}
\includegraphics[width=5cm,height=3cm]{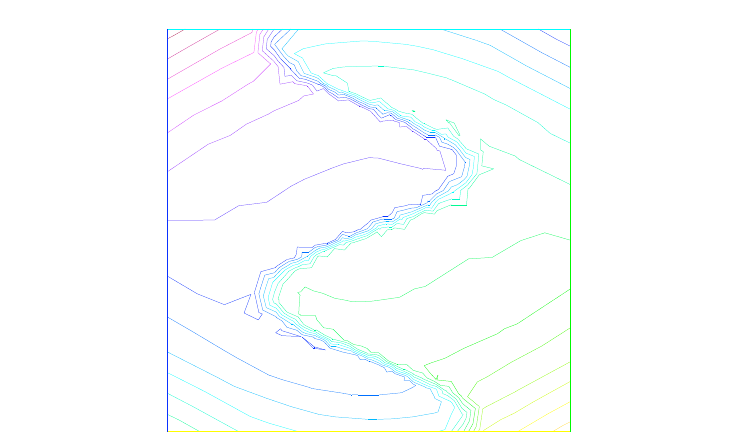}
\hspace{-20mm}
\includegraphics[width=5cm,height=3cm]{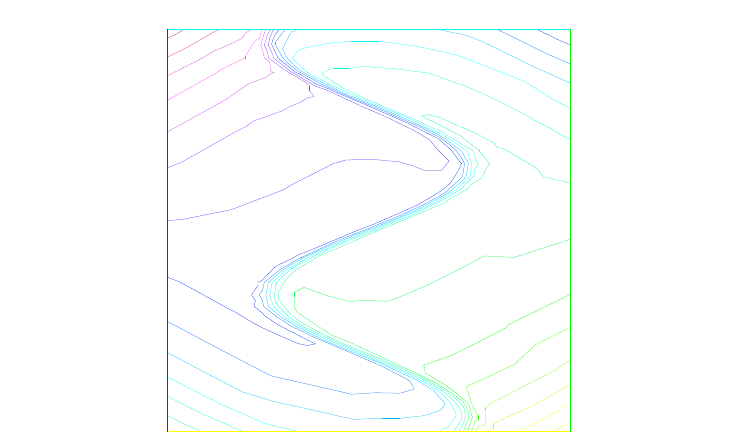}
\hspace{-30mm}
\vspace{-3mm}
\caption{Interpolation of \iref{syntheticF} with $\P_1$ elements on a uniform, isotropic or anisotropic mesh of cardinality $500$.}
\end{figure}

%\bibliographystyle{plain}
%\bibliography{myRef}

\end{document}